\documentclass[11pt,a4paper]{amsart}
\usepackage[all]{xy}
\usepackage{amssymb}

\numberwithin{equation}{section}

\newcommand{\ie}{{\em i.e.}\ }
\newcommand{\cf}{{\em cf.}\ }

\newcommand{\ko}{\: , \;}

\newtheorem{theorem}{Theorem}

\newtheorem{lemma}{Lemma}
\newtheorem{proposition}{Proposition}
\newtheorem{corollary}{Corollary}

\newcommand{\opname}[1]{\operatorname{\mathsf{#1}}}

\renewcommand{\mod}{\opname{mod}\nolimits}

\newcommand{\Gr}{\opname{Gr}\nolimits}
\newcommand{\dimv}{\underline{\dim}\,}

\newcommand{\ac}{\mathcal{A}}

\newcommand{\cok}{\opname{cok}\nolimits}
\newcommand{\im}{\opname{im}\nolimits}
\renewcommand{\ker}{\opname{ker}\nolimits}

\newcommand{\Z}{\mathbb{Z}}
\newcommand{\N}{\mathbb{N}}
\newcommand{\Q}{\mathbb{Q}}

\renewcommand{\P}{\mathbb{P}}

\newcommand{\la}{\leftarrow}
\newcommand{\iso}{\stackrel{_\sim}{\rightarrow}}

\newcommand{\Hom}{\opname{Hom}}
\newcommand{\go}{\opname{G_0}}
\newcommand{\gr}{\opname{Gr}}

\newcommand{\Ext}{\opname{Ext}}

\newcommand{\Ind}{\opname{Ind}}
\newcommand{\End}{\opname{End}}
\newcommand{\rad}{\opname{rad}}
\newcommand{\GL}{\opname{GL}}

\newcommand{\ca}{{\mathcal A}}

\newcommand{\cc}{{\mathcal C}}
\newcommand{\cd}{{\mathcal D}}

\newcommand{\cF}{{\mathcal F}}

\newcommand{\ch}{{\mathcal H}}

\newcommand{\cm}{{\mathcal M}}

\newcommand{\cp}{{\mathcal P}}

\newcommand{\cu}{{\mathcal U}}
\newcommand{\cv}{{\mathcal V}}

\newcommand{\db}{{\mathcal D}^b}

\newcommand{\hmod}{H\opname{-mod}}
\newcommand{\den}{\opname{den}\nolimits}

\newcommand{\bt}{\bullet}
\newcommand{\sq}{\square}
\newcommand{\point}{.}
\newcommand{\dmd}{\diamondsuit}
\newcommand{\sqdmd}{\square\!\!\!\!\!\raisebox{1pt}{$\diamondsuit$}}
\newcommand{\sqbt}{\square\!\!\!\!\bt}
\newcommand{\sqpt}{\square\!\!\!\raisebox{3pt}{\point}}

\begin{document}

\title{From triangulated categories to cluster algebras II\\}
\author{Philippe Caldero}
\address{Institut Camille Jordan, Universit\'e Claude Bernard Lyon I,
69622 Villeurbanne Cedex, France}
\email{caldero@math.univ-lyon1.fr}
\author{Bernhard Keller}
\address{UFR de Math\'ematiques\\
   Institut de Math\'ematiques \\
   Case 7012\\
   Universit\'e Paris 7 - Denis Diderot\\
   2, place Jussieu\\
   75251 Paris Cedex 05\\
   France }
\email{keller@math.jussieu.fr}

\begin{abstract}
In the acyclic case, we establish a one-to-one correspondence
between the tilting objects of the cluster category and the
clusters of the associated cluster algebra. This correspondence
enables us to solve conjectures on cluster algebras. We prove a
multiplicativity theorem, a denominator theorem, and some conjectures on
properties of the mutation graph. As in the previous article, the
proofs rely on the Calabi-Yau property of the cluster category.
\end{abstract}
\date{version du 12/10/2005, modifi\'ee le 12/04/2006}

\maketitle
\section{Introduction}
Cluster algebras are commutative algebras, introduced in
\cite{fominzelevinsky1} by S.~Fomin and A.~Zelevinsky. Originally,
they were constructed to obtain a better understanding of the
positivity and multiplicativity properties of Lusztig's dual
(semi)canonical basis of the algebra of coordinate functions on
homogeneous spaces. Cluster algebras are generated by the
so-called {\it cluster variables} gathered into sets of fixed
cardinality called {\it clusters}. In the framework of the present
paper, the cluster variables are obtained by a recursive process
from an antisymmetric square matrix $B$.

Denote by $Q$ the quiver associated to the matrix $B$. Assume that $Q$
is connected. A theorem of Fomin and Zelevinsky asserts that the
number of cluster variables of the corresponding cluster algebra
$\ca_Q$ is finite if and only if $Q$ is mutation-equivalent to a
quiver whose underlying graph is a simply laced Dynkin diagram. In
this case, it is known that the combinatorics of the clusters are
governed by the generalized associahedron.

Let $Q$ be any finite quiver without oriented cycles and let $k$
be an algebraically closed field. The cluster category $\cc=\cc_Q$
was introduced in \cite{CCS} for type A$_n$  and in \cite{BMRRT}
in the general case. This construction was motivated by the
combinatorial similarities of $\cc_Q$ with the cluster algebra
$\ca_Q$. The cluster category is the category of orbits under an
autoequivalence of the bounded derived category $\db$ of the
category of finite dimensional $kQ$-modules. By \cite{keller}, the
category $\cc_Q$ is a triangulated category. Let us denote its
shift functor by $S$ and write $\Ext^1_\cc(M,N)$ for $\Hom_\cc(M,SN)$
for any objects $M,N$ of $\cc$. By construction, the cluster category is
Calabi-Yau of CY-dimension $2$; in other terms, the functor
$\Ext^1$ is symmetric in the following sense:
 $$\Ext_{\cc}^1(M,N)\simeq D\Ext_{\cc}^1(N,M).$$

In a series of articles \cite{BMRRT}, \cite{BMR}, \cite{BMR2}, the
authors study the tilting theory of the cluster category. More
precisely, they describe the combinatorics of the cluster tilting
objects of the category $\cc$, \ie the objects without
self-extensions and with a maximal number of non-isomorphic
indecomposable summands. In \cite{BMR2}, the authors define a map
$\beta$ between the set of clusters of $\ac_Q$ and the set of
tilting objects of the category $\cc_Q$. A natural question
arises: does  $\beta$ provide a one-to-one correspondence between
both sets?

In the articles \cite{caldchap} and \cite{caldkell}, it is proved
that in the finite case, \ie the Dynkin case, the cluster algebra
can be recovered from the corresponding cluster category as the
so-called {\it exceptional Hall algebra} of the cluster category.
More precisely, in \cite{caldchap}, the authors give an explicit
correspondence $M\mapsto X_M$ between indecomposable objects of
$\cc_Q$ and cluster variables of $\ca_Q$. In \cite{caldkell}, we
provide a multiplication rule for the algebra $\ca_Q$ in terms of
the triangulated category $\cc_Q$.

An ingenious application of the methods of \cite{caldkell} can be
found in \cite{GLS}, where the authors give a multiplication
formula for elements of Lusztig's dual semicanonical basis. Here,
the cluster category is replaced by the category of
finite-dimensional modules over the preprojective algebra and the
r\^ole of the cluster algebra is played by the coordinate algebra
of the maximal unipotent subgroup in the corresponding semisimple
algebraic group.

The aim of the present article is  to generalize some of the
results of \cite{caldchap}, \cite{caldkell} to the case where $Q$
is any finite quiver without oriented cycles. Building on the
important results obtained in \cite{BMR2} we strengthen here the
connections between the cluster category and the cluster algebra
by giving an explicit expression for the correspondence $\beta$
and proving that $\beta$ is one-to-one. The key ingredient of the
proof is a natural analogue of the map $M\mapsto X_M$ of
\cite{caldchap}. With the help of a multiplicativity result, we show
that $M\mapsto X_M$ defines a bijection between the indecomposable
objects without self-extensions of $\cc_Q$ and the cluster
variables of $\ac_Q$.

This correspondence between cluster algebras and cluster
categories gives positive answers to some of the conjectures which
S.~Fomin and A.~Zelevinsky formulated in \cite{FZ-conf}. We prove connectedness
properties of some mutation graphs, \cf
section~\ref{ss:connectedness}. 
As a byproduct, we obtain a
cluster-categorical interpretation of the passage to a submatrix
of the exchange matrix. This strengthens a key result of
\cite{BMR2} and may be of independent interest.

Another consequence of the bijectivity of $\beta$ is that
each seed is determined by its cluster. As we have learned
recently, this result is obtained independently in \cite{BMRT}.

The paper is organized as follows: In the first part, we recall
well-known facts on the cluster category. For any object $M$ of
the cluster category, we define the Laurent polynomial $X_M$ 
as in \cite{caldchap}. 
With the techniques of \cite{caldkell}, we prove
an `exchange relation' for the $X_M$. To be more precise, we prove
that if $M$ and $N$ are indecomposable objects of the category
$\cc=\cc_Q$ such that $\Ext_{\cc}^1(M,N)=k$, then
\[
X_MX_N=X_B+X_{B'},
\]
where $B$ and $B'$ are the unique objects (up to isomorphism) such
that there exist non split triangles
\[
N\rightarrow B\rightarrow M\rightarrow SN, \;\;
M\rightarrow B'\rightarrow N\rightarrow SM.
\]
This formula is an analogue of the `exchange relation' between cluster
variables. With the help of a comparison theorem of \cite{BMR2}, we
prove by induction that for any indecomposable exceptional object
$M$, $X_M$ is a cluster variable and that its (monomial) denominator
is given by the dimension vector $\dimv(M)$.  From this denominator
property, we deduce that the map $M\mapsto X_M$ is injective when
restricted to the set of indecomposable objects of $\cc_Q$ without
self-extensions. The connectedness of the tilting graph
proved in \cite{BMR2} then implies
that the map $M\mapsto X_M$ is a one-to-one correspondence between the
set of tilting objects of $\cc_Q$ and the set of clusters of $\ac_Q$.

We then deduce some applications of this correspondence to conjectures
of \cite{FZ-conf}.

 {\bf Acknowledgements:} The first author is indebted to Thomas
Br\"ustle and  Ralf Schiffler for useful
conversations. He also wishes to thank Andrei Zelevinsky for his kind
hospitality and for pointing out to him the conjectures of
\cite{FZ-conf}. The authors thank Andrew Hubery for pointing
out a gap in a previous version of this article.

\newpage
\section{The cluster category and the cluster variable formula}

\begin{subsection}{}\label{H-mod} Let $H$ be a finite dimensional
hereditary algebra over an algebraically closed field $k$. We
denote by $\hmod$ the category of finitely generated $H$-modules.
We choose representatives $S_i$, $1\leq i\leq n$, of the
isoclasses of the simple $H$-modules and denote by $I_i$ the
injective hull and by $P_i$ the projective cover of $S_i$.

The Grothendieck group of $\hmod$ is the group $\go(\hmod)$
generated by the isoclasses of  modules in $\hmod$ and subject to
the relations $X=M+N$ obtained from exact sequences $0\rightarrow
M\rightarrow X\rightarrow N\rightarrow 0$ in $\hmod$. We denote by
$[M]$ the class of a module $M$ in $\go(\hmod)$. We put
$\alpha_i=[S_i]$. The Grothendieck group is free abelian on the
$\alpha_i$. The dimension vector $\dimv(M)$ of a module $M$ is by
definition the vector of the coordinates of $[M]$ in this basis.

We define the Euler form by $<M,N>=\dim\Hom(M,N)-\dim\Ext^1(M,N)$,
for any $M$, $N$ in $\hmod$. Since $H$ is hereditary, this form is
well-defined on the Grothendieck group.

Let $\tau$ be the Auslander-Reiten functor of $\hmod$. This
functor verifies the Auslander-Reiten formula:
\[ D\Hom(N,\tau M)=\Ext^1(M,N),\]
where $D$ is the functor $\Hom_k(?,k)$.
\end{subsection}

\begin{subsection}{}
For any $H$-module $M$, and any $e$ in $\go(\hmod)$, we denote by
$\gr_e(M)$ the Grassmannian of submodules of $M$ with dimension
vector $e$:
\[\gr_e(M)=\{N,\,N\in\hmod,\,N\subset M,\, \dimv(N)=e\}.\]
It is a closed subvariety of the classical Grassmannian of the
vector space $M$. Let $\chi_c$ be the Euler-Poincar\'e
characteristic of the etale cohomology with proper support defined
by
$$\chi_c(X) = \sum_{i=0}^\infty (-1)^i \dim H^i_c(X, \overline\Q_l).$$

Let $\Q[x_i^{\pm 1}, 1\leq i\leq n]$ be the $\Q$-algebra of
Laurent polynomials in the variables $x_i$. As in
\cite{caldchap}, for any module $M$, we set
\[X_M=\sum_e\chi_c(\gr_e(M))\prod_i
x_i^{-<e,\alpha_i>-<\alpha_i,m-e>}\in \Q[x_i^{\pm 1}, 1\leq i\leq
n],\] where $m:=\dimv(M)$. Note that, as $M$ is finite
dimensional, there only exists a finite number of non zero terms
in this sum. Remark that $X_M$ only depends on the isoclass of the
module $M$. As in \cite{caldchap} one shows that
\[
\chi_c(\gr_g(M\oplus
N))=\sum_{e+f=g}\chi_c(\gr_e(M))\chi_c(\gr_f(N)).
\]
Hence, the bilinearity of the Euler form implies that
\[X_{M\oplus N}=X_M X_N.\]
\end{subsection}

\begin{subsection}{}
As $H$ is hereditary and finite dimensional, there exists a finite
quiver $Q$ without oriented cycles such that $H$ is Morita equivalent
to the path algebra $kQ$ of $Q$. Let $Q_0$ be the set of vertices and
$Q_1$ the set arrows of $Q$. Let $n$ be the number of vertices of
$Q$.

 The bounded derived category $\db=\db(H)$ of $\hmod$ is a triangulated
category. We denote its shift functor $M\mapsto M[1]$ by  $S$.
The category $\db$ is a Krull-Schmidt category and, up to
canonical triangle equivalence, it only depends on the underlying
graph of $Q$, see \cite{happel}. We identify the category $\hmod$ with
the full subcategory of $\db$ formed by the complexes whose homology
is concentrated in degree $0$.  We simply call `modules' the objects
in this subcategory. The indecomposable objects of $\db$ are the
shifts $S^iM$, $i\in\Z$, of the indecomposable objects of $\hmod$. We
still denote by $\tau$ the AR-functor of $\db$; it is known that
$\tau$ is an autoequivalence characterized by the Auslander-Reiten
formula.

Let $F$ be the autoequivalence $\tau^{-1}S$ of $\db$. The AR-formula
implies that
\[\Ext_{\db}^1(M,N)=\Hom_{\db}(M,SN)=D\Ext_{\db}^1(FN,M),\]
for any objects $M$, $N$ of $\db$. Let $\cc=\cc(H)$ be the orbit
category $\db/F$: the objects of $\cc$ are the objects of $\db$ and
the morphisms of $\cc$ are given by
\[\Hom_{\cc}(M,N)=\oplus_{i\in \Z}\Hom_{\db}(M,F^iN).\]
The category $\cc$ is the so-called cluster category, introduced and studied
in depth in \cite{BMRRT}.  Let $\pi$ be the canonical functor from $\db$ to
$\cc$. We will often omit the functor $\pi$ from the notations. 
Statements (i) and (ii) of the following theorem were proved
in \cite{keller}, statements (iii) and (iv) in \cite{BMRRT}:
\begin{theorem}\label{basic}
\begin{itemize}
  \item[(i)] The category $\cc$ is triangulated and
  \item[(ii)] the functor $\pi$ : $\cd\rightarrow\cc$ is
   a triangle functor.
  \item[(iii)] The category $\cc$ is a Krull-Schmidt category and
  \item[(iv)]  For any indecomposable object without self-extensions
               $M$ of $\cc$, we have $\End_{\cc}(M)=k$.
\end{itemize}
\end{theorem}
The shift functor of the triangulated category $\cc$ will still be
denoted by $S$.  For any objects $M$, $N$ of $\cc$, the formulas
above imply that there exists an (almost canonical) duality
\[\phi\,:\,\Ext_{\cc}^1(M,N)\times\Ext_{\cc}^1(N,M)\rightarrow k.\]
The set of indecomposable objects of $\cc$ is given by
\[\Ind(\cc)=\Ind(\hmod)\coprod\{SP_i,\,1\leq i\leq n\}.\]
Note that $SP_i=S^{-1}\tau^{-1}SP_i=S^{-1}I_i.$

We extend the definition of $X_M$ to any object $M$ of the
category $\cc$ by setting $X_{SP_i}=x_i$, $1\leq i\leq n$, and
requiring $X_{M\oplus N}=X_MX_N$ for all objects $M,N$ of $\cc$.

The AR-formula and the fact that $\tau$ passes to the Grothendieck
group of the derived category of $\hmod$ allow us to rewrite $X_M$
for a module $M$ as
\begin{equation}\label{ccformula}
X_M = \sum_{e} \chi_c(\Gr_e(M)) \, x^{\tau(e) - \dimv M + e} \ko
\end{equation}
where we have set
\[
x^v = \prod_{i=1}^n x_i^{<\dimv S_i, v>},
\]
for any $v$ in $\Z^n$. Remark that this notation gives
\[X_{SP_i}=x^{\dimv I_i}.\]

\end{subsection}
\begin{subsection}{}
Each object $M$ of $\cc$ can be uniquely decomposed in
the following way:
$$M=M_0\oplus SP_M=M_0\oplus S^{-1}I_M,$$ where $M_0$ is the image
under $\pi$ of a module in $\db$, and where $P_M$, respectively $I_M$, is a
uniquely determined projective, respectively injective, module. We will say
that an object $M$ of $\cc$ is a {\em module} if $M=M_0$, and that $M$
is the {\em shift of a projective module} if $M=SP_M$.\par

From \cite{BMRRT}, we recall the
\begin{proposition}\label{symmetric}
For any indecomposable modules $M$ and $N$ in $\cc$, we have
\[\Ext_{\cc}^1(M,N)=\Ext_{H}^1(M,N)\oplus D\Ext^1_H(N,M).\]
\end{proposition}
The module $M_0$ can be recovered using the functor
$$H^0=\Hom_{\cc}(H_{H},?)\,:\,\cc\rightarrow \hmod.$$
Indeed, we have
\[
H^0(M)=H^0(M_0)\oplus H^0(SP_M)=\Hom_{\hmod}(H_{H},M_0)\oplus
\Hom_{\cc}(\oplus_i P_i,SP_M)=M_0 \ko
\]
as the last factor is zero. The functor $H^0$ is a homological
functor, \ie it maps triangles in $\cc$ to long exact sequences of
$H$-modules.
\end{subsection}

\section{A multiplication formula}
\begin{subsection}{}
The aim of the section is to prove the following theorem:
\begin{theorem}\label{exchange}
Let $M$ and $N$ be indecomposable objects of the category $\cc$ such
that $\Ext_{\cc}^1(M,N)$ is one-dimensional. Then we have
\[X_MX_N=X_B+X_{B'},\]
where $B$ and $B'$ are the unique objects (up to isomorphism) such
that there exist non split triangles
\[N\rightarrow B\rightarrow M\rightarrow SN, \;\;
M\rightarrow B'\rightarrow N\rightarrow SM.\]
\end{theorem}
Note that when $H$ is the path algebra of a Dynkin quiver,
the theorem is a particular case of the cluster multiplication
formula of \cite{caldkell}. Actually, we will see that the method of
\cite{caldkell} generalizes nicely to the framework of the theorem.

Thanks to the hypotheses of the theorem and the symmetry of $\Ext^1$,
we just need to consider the two following cases

1. $N=SP_i$ for an $i\in Q_0$ and $M$ is an indecomposable module.

2. $M$ and $N$ are indecomposable modules.

\noindent
Indeed, the isomorphisms $M=SP_j$ and $N=SP_i$ would imply
\[
\Ext_{\cc}^1(M,N)=\Ext_{\cc}^1(P_j,P_i)=0.
\]

\end{subsection}
\begin{subsection}{}
We now prove the theorem in the first case. Suppose $N=SP_i$, and let
$M$ be an indecomposable module such that
$\Ext_{\cc}^1(SP_i,M)=k$. Using theorem \ref{basic} and the
AR-formula, we obtain
\[
\dimv(M)_i=\dim\Hom_H(P_i,M)=\dim\Hom_{\cc}(P_i,M)=\dim\Ext_{\cc}^1(SP_i,M)=1.
\]
Hence, up to a multiplicative scalar, there exists a unique non zero
morphism $\zeta$ : $M\rightarrow I_i$ and a non zero morphism
$\zeta': P_i\rightarrow M$.
\begin{lemma}
Let $M'$ be a submodule of $M$. Then either $M'\subset\ker\zeta$ or
$\im\zeta'\subset M'$.
\end{lemma}
\begin{proof}
By the formula above, the space $M'_i$ is of dimension $0$ or $1$. We claim
that

1. $\dim(M')=0$ if and only if $M'\subset\ker\zeta$,

2. $\dim(M')=1$ if and only if $\im\zeta'\subset M'$.

The lemma follows from the claim. Let's prove part 1. The second part
is similar and left to the reader.  The module $\im\zeta$ is non zero
and so it contains the simple socle $S_i$ of $I_i$. Hence,
$\dim(\ker\zeta)_i=0$, which gives the `if' part.
Conversely, if $\dim(M')_i=0$, then $\zeta(M')\cap S_i=0$, hence
$\zeta(M')=0$ as $S_i$ is the socle of $I_i$.
\end{proof}

Applying the functor $H^0$ the  non split triangle
\begin{equation}\label{triangle}
\xymatrix{SP_i\ar[r]^\iota &B\ar[r]^\pi &M\ar[r]^-\zeta &S^2P_i=I_i}
\end{equation}
we obtain a long exact sequences of $H$-modules
\[
\xymatrix{0\ar[r]&H^0B \ar[r]^{H^0\pi} & M \ar[r]^{H^0\zeta} &
I_i \ar[r]^-{H^1\iota}&H^0\tau B\ar[r]^{H^1\pi}&H^0\tau M} \ko
\]
Now, $H^0\tau B=\tau H^0B\oplus I_B$, and the first factor is non
injective. As the quotient of an injective module is still injective,
we have $\im(H^1\iota)\subset I_B$. Moreover, as $H^0\tau M$ is non
injective, we have $I_B\subset\ker(H^1\pi)$. Hence, we have equality
and so the following exact sequence holds
\begin{equation}\label{exact}
\xymatrix{0\ar[r]&H^0B \ar[r]^{H^0\pi} & M \ar[r]^{H^0\zeta} & I_i
\ar[r]^{H^1\iota}&I_B\ar[r]&0}.
\end{equation}
Note that the morphism $H^0\zeta = \zeta$ is non zero.

In the same way, applying the functor $H^0$, the  non split triangle
\begin{equation}\label{triangle'}
\xymatrix{P_i\ar[r]^{\zeta'} &M\ar[r]^{\pi'} &B'\ar[r]^{\iota'} &SP_i,}
\end{equation}
we obtain
\begin{equation}\label{exact'}
\xymatrix{0\ar[r]&P_{B'} \ar[r] & P_i \ar[r]^{H^0\zeta'} &
M \ar[r]^{H^0\pi'}&H^0B'\ar[r]&0}.
\end{equation}
Note that the morphism $H^0\zeta'=\zeta'$ is non zero.

Now, the lemma implies that for any submodule $M'$ of $M$, $M'$ is a
either submodule of $\im H^0\pi$ or contains $\ker H^0\pi'$. Hence,
there is a natural bijection between $\gr_eM$ and
$\gr_e(H^0B)\coprod\gr_{e-k}(H^0B')$, where
\begin{equation}\label{dimvker}
k:=\dimv\ker H^0\pi'=\dimv P_i-\dimv P_{B'}.
\end{equation}
We want to prove the multiplication formula, which in this case is
\[x^{\dimv I_i}\sum_e\chi_c(\gr_eM)x^{\tau(e)-\dimv M+e}=\hskip 5cm\]
\[x^{\dimv I_B}\sum_e\chi_c(\gr_e H^0B)x^{\tau(e)-\dimv H^0B+e}+
x^{\dimv I_{B'}}\sum_e\chi_c(\gr_eH^0B')x^{\tau(e)-\dimv H^0B'+e}.\]
So, it remains to prove that
\[\dimv I_i+\tau(e)-\dimv M+e=\dimv I_B+\tau(e)- \dimv H^0B+e,\]
and
\[\dimv I_i+\tau(e)-\dimv M+e=\dimv I_{B'}+\tau(e-k)- \dimv H^0B'+e-k.\]
The first formula is a direct consequence of \ref{exact}. The second
one comes from \ref{exact'}, \ref{dimvker} and the formula $\tau(\dimv
P_j)=-\dimv I_j$.
\end{subsection}
\begin{subsection}{}
This subsection and the following one are devoted to the proof of the
theorem in the second case.  In order to simplify notations, we will
write $(X,Y)$ for $\Hom_{\cc}(X,Y)$.

Let $M$ and $N$ be two indecomposable modules such that
$\Ext^1_{\cc}(N,M)=k$. By proposition \ref{symmetric}, we can suppose
that $\Ext_H^1(N,M)=k$ and $\Ext_H^1(M,N)=0$. In this case, by theorem
\ref{basic}, there exists (up to isomorphism) a unique non split short
exact sequence of $H$-modules

\[\xymatrix{0\ar[r]&M \ar[r]^i & B_+ \ar[r]^p & N \ar[r]&0},\]
and two triangles in $\cc$
\[\xymatrix{M \ar[r]^i & B_+ \ar[r]^p & N \ar[r]&SM},\]
\[\xymatrix{N \ar[r]^{i'} & B_- \ar[r]^{p'} & M \ar[r]&SN}.\]
Note that $B_+$ is a `module' of $\cc$ but $B_-$ is just an object;
they both are uniquely determined up to isomorphism. We want to prove
the formula
\[X_{B_+}+X_{B_-}=X_MX_N,\]
and the idea is first to construct a morphism $\Psi$ between $\gr
B_+\coprod\gr H^0B_-$ and $\gr M\times\gr N$. For any submodule $B_+'$
of $B_+$, set $\Psi(B_+')=(i^{-1}B_+',pB_+')$, and for any submodule
$B_-'$ of $B_-$, set $\Psi(B_-')=((H^0p')B_-',(H^0i')^{-1}B_-')$.  As
a first step, we want to prove the proposition
\begin{proposition}\label{surjective}
The variety $\gr M\times\gr N$ is the disjoint union of $\Psi(\gr
B_+)$ and $\Psi(\gr H^0B_-)$. Moreover, the fibers of $\Psi$ are affine
spaces.
\end{proposition}
This proposition will be proved at the end of this subsection.

Given a submodule $M'$ of $M$, a submodule $N'$ of $N$, and the
corresponding embeddings $i_M$ and $i_N$, we have a diagram
\[
\xymatrix{ S^{-1}M \ar[r]^-{\varepsilon'} & N \ar[r]^{\varepsilon} & SM
              \\ S^{-1}M' \ar[u]^{S^{-1}i_M} & N' \ar[u]^{i_N} & SM'
              \ar[u]^{Si_M}}
\]
and two complexes
\[
\xymatrix{
(S^{-1}M,N') \ar[r]^-{\alpha'} & (S^{-1}M,N)\oplus (S^{-1}M',N')
\ar[r]^-{\beta'} & (S^{-1}M',N) \\
(N',SM) & (N,SM)\oplus (N',SM') \ar[l]_-{\alpha} & (N,SM')
\ar[l]_-{\beta}}
\]
where
\[
\alpha' = \left[ \begin{array}{c} (i_{N'})_* \\ (S^{-1} i_{M'})^*
\end{array} \right] \ko
\beta' = \left[ (S^{-1} i_{M'})^*, -(i_N')_*  \right] \ko
\alpha = \left[ (i_{N'})_*, (S i_{M'})^* \right] \ko
\beta = \left[ \begin{array}{c} (S i_{M'})^* \\ -(i_{N'})_* \end{array}
\right] .
\]
The two sequences are dual to each other via the canonical duality $\phi$.

The following proposition is straightforward by using basic properties of triangulated categories.
\begin{proposition}\label{commute}
The following conditions are equivalent:
\begin{itemize}
\item[(i)] There exists a submodule $B_+' \subset B_+$ such that the diagram
  \[
\xymatrix{ 0 \ar[r] & M \ar[r] & B_+ \ar[r] & N
\ar[r] & 0 \\
           0 \ar[r] & M' \ar[u]\ar[r] & B_+' \ar[u]\ar[r] &
           N' \ar[u]\ar[r] & 0 .
}
\]

  commutes.
  \item[(ii)] There exists a morphism $\eta$ : $N'\rightarrow SM'$
such that the square
  \[
\xymatrix{   N \ar[r]^{\varepsilon}         & SM \\
             N' \ar[u]^{i_N}\ar[r]^\eta  & SM' \ar[u]^{Si_M}}
\]
commutes.
  \item[(iii)] The composed morphism
  \[\ker\alpha\hookrightarrow (N',SM')\oplus (N,SM)\rightarrow (N,SM)\]
  is non zero.
  \item[(iv)]  The composed morphism
  \[\cok\alpha'\leftarrow (S^{-1}M',N')\oplus (S^{-1}M,SN)\hookleftarrow
   (S^{-1}M,N)\]
  is non zero.
\end{itemize}

\end{proposition}
The following proposition sheds light on the situation when the
conditions of proposition \ref{commute} do not hold.
\begin{proposition}\label{ifnot}
The following conditions are equivalent:
\begin{itemize}
  \item[(i)] The composition
  \[\cok\alpha'\leftarrow (S^{-1}M',N'')\oplus (S^{-1}M,SN)\hookleftarrow
   (S^{-1}M,N)\]
   vanishes, \ie $(S^{-1}M,N)$ is contained in the image of $\alpha'$.
  \item[(ii)] There exist a submodule $B_-'\hookrightarrow H^0B_-$ and
  a commutative diagram
  \[
\xymatrix{ N \ar[r]^{H^0i'} & H^0B' \ar[r]^{H^0p'} & M&\\N'
\ar[u]\ar[r] & B_-' \ar[u]\ar[r] &M' \ar[u]\ar[r]&0
\\0\ar[u]&0\ar[u]&0\ar[u]&}\]
where $M':=H^0p'(B_-')$ and $N'=(H^0i')^{-1}(B_-')$.
\end{itemize}

\end{proposition}
\begin{proof} Let us show that (i) implies (ii).
By the assumption, we can find a commutative square
\[
\xymatrix{S^{-1}M\ar[rr]^{\varepsilon'}
\ar[drr]^{f} & &
    N\\S^{-1}M' \ar[u]^{S^{-1}i_M} \ar[rr]_0 & &
    N' \ar[u]^{i_N}}
\]
We complete it to a morphism between triangles:
\[
\xymatrix{S^{-1}M \ar[rr]^{\varepsilon'}
\ar[drr]^{f} & &
    N \ar[r]^{i'} & B_- \ar[r]^{p'} & M \\S^{-1}M' \ar[u] \ar[rr]_0 & &
    N' \ar[u]^{i_N} \ar[r] & N'\oplus M' \ar[u]\ar[r] & M' \ar[u]^{i_M}}
\]
We take the homology:
\[
\xymatrix{H^0(S^{-1}M) \ar[rr]^{H^0(\varepsilon')}
\ar[drr]^{H^0(f)} & &
    N \ar[r]^{H^0 i'} & H^0 B_- \ar[r]^{H^0 p'} & M&\\0\ar[rr]&&N' \ar[u]^{i_N} \ar[r] &
    N'\oplus M' \ar[u] \ar[r] & M' \ar[u]^{i_M}\ar[r]&0}
\]
We take $B_-'$ as the image of $N'\oplus M'\rightarrow H^0B_-$. Let us
show that $N'\subset N$ is $H^0(i')^{-1}(B_-')$. Indeed, clearly the
image of $N'$ is contained in $B_-'$. Conversely, if we have $x\in N$
whose image lies in $B_-'$, then the image is the image of $(x',y')$
in $N'\oplus M'$, and the image of $x\in N$ under $N\rightarrow
H^0B_-\rightarrow M $ vanishes. So, the image of $y'$ in $M$
vanishes. But $M'\rightarrow M$ is mono. So $y'$ vanishes and we get
$x'$ in $N'$ such that $x$ in $N$ and $x'$ have the same image in
$H^0B_-$. Then $x=x'+(H^0(\varepsilon')(z)$ for some $z$ in
$S^{-1}M$. But $H^0\varepsilon'=(H^0i_N)\circ(H^0f)$. So
$(H^0\varepsilon')(z)$ lies in fact in $N'\subset N$ and $x$ lies in
$N'$.

Let us show that $M'$ is the image of $B_-'$. Clearly, the image of
$B_-'$ is contained in $M'$. Conversely, if $x'\in M'$, we consider
the image $y$ in $B_-'$ of $(0,x')\in N'\oplus M'$. Then clearly, the
image of $y$ is $x'$.

Let us prove that (ii) implies (i).
The hypothesis yields the following diagram
\[
\xymatrix{H^0(S^{-1}M) \ar[rr]^{H^0(\varepsilon')}
\ar@{.>}[drr] & &
    N \ar[r]^{H^0 i'} & H^0B_- \ar[r]^{H^0 p'} & M \ar[r]^-{H^0(S\varepsilon')} &   H^0(SN)\\
H^0(S^{-1}M') \ar[u]  & &
    N' \ar[u] \ar[r] & B_-' \ar[u]\ar[r] & M' \ar[u] \ar[r]_-0   & H^0(SN').
    \ar[u]\\
    &&0\ar[u]&0\ar[u]&0\ar[u]&0\ar[u]}
\]

As the composition $H^0(S^{-1}M)\rightarrow N\rightarrow H^0B_-$
vanishes, the image of $H^0(S^{-1}M)$ is contained in $N'$, which is
the inverse image of $B_-'$.

As the composition $H^0(B_-)\rightarrow M\rightarrow H^0(SN)$
vanishes, $M'$ is contained in the kernel of $M\rightarrow
H^0(SN)$. We know that $M$ is not injective, so, $S^{-1}M=\tau^{-1}M$
is still a module.

Moreover, we have
\[D\Ext^1(N,M)=\Hom_H(\tau^{-1}M,N)=\Hom_H(M,\tau N)=k.\]

We obtain the commutative diagrams

\[ \xymatrix{S^{-1}M \ar[rr]^{\varepsilon'}
\ar@{.>}[drr]_f & &
    N\\  & &N'\ar[u]_{i_N}}\hskip 1cm
 \xymatrix{M \ar[rr]^{S\varepsilon'}
 & &SN\\  M'\ar[u]^{i_M}\ar@{.>}[urr]_0& &} \]

The module $M'$ has no injective direct summand, because $M$ is
indecomposable and non injective. So, $S^{-1}M'$ is still a
module. Consider

\[
\xymatrix{S^{-1}M \ar[rr]^{\varepsilon'} \ar[drr]^{f} & & N \\S^{-1}M
\ar[u]_{S^{-1}i_M} \ar[rr]_0 & & N'\ar[u]^{i_N}}.
\]
We have $i_N\circ f\circ S^{-1}i_M=\varepsilon'\circ S^{-1}i_M=0$. As
$i_N$ is injective, this gives $f\circ S^{-1}M=0$, which implies (i).
\end{proof}

Propositions \ref{commute} and \ref{ifnot} imply the first part of
proposition \ref{surjective}. The second part is a well-known
fact, \cf lemma~3.8 of \cite{caldchap}.
\end{subsection}
\begin{subsection}{}

We want to prove the multiplication formula for the second. It reads
as follows:
\[
\sum_e\chi_c(\gr_eM)x^{\tau(e)-\dimv M+e}\sum_f\chi_c(\gr_fN)x^{\tau(f)-\dimv N+f}=\hskip 2cm\]
\[\sum_g\chi_c(\gr_g H^0B_+)x^{\tau(g)-\dimv B_++g}+x^{\dimv I_{B_-}}\sum_g\chi_c(\gr_gH^0B_-)x^{\tau(g)-\dimv H^0B_-+g}.
\]
By combining proposition~\ref{surjective} with proposition~3.6 of
\cite{caldchap}, we can compare Euler characteristics on both sides of
the equality. What we need to prove now is
\begin{equation}\label{first}
\tau(e)-\dimv M+e+\tau(f)-\dimv N+f=\tau(g)-\dimv B_++g,
\end{equation}
with $e=\dimv M'$, $f=\dimv N'$, $g=\dimv B_+'$, in the setting of
proposition~\ref{commute} (i), and then
\begin{equation}\label{second}
\tau(e)-\dimv M+e+\tau(f)-\dimv N+f=\dimv I_{B_-}+\tau(g)-\dimv H^0B_-+g,
\end{equation}
with $e=\dimv M'$, $f=\dimv N'$, $g=\dimv B_-'$, in the setting of
proposition~\ref{ifnot} (ii).

The  formula~\ref{first} is clear since $g=e+f$ in this case.

In order to prove the second formula, we need to complete the diagram
of proposition \ref{ifnot} by adding kernels and cokernels

 \[
\xymatrix{     &     & N/N' \ar[r]        & (H^0 B_-)/B_-' \ar[r]  & M/M' \ar[r]&C\ar[r]&0  \\
&H^0(S^{-1}M) \ar[r] & N \ar[u]\ar[r]^{H^0(i)}     & H^0 B_- \ar[u]\ar[r]^{H^0(p)} & M
\ar[u] \ar[r] & H^0(SN)& \\
   0\ar[r] &K\ar[r]                & N' \ar[u]\ar[r]    & B_-' \ar[u]\ar[r] & M' \ar[u] &&
}
\]

With the notation above, the diagram implies the equalities

\[
\tau(e)+\tau(f)=\tau(g)+\tau(\dimv K),\;\;(\dimv M-e)+(\dimv N -f)=
(\dimv H^0B_- -g)+\dimv C.
\]
So, in order to prove formula \ref{second}, we remains  to show that
\begin{equation}\label{wanted}
\tau\dimv  K-\dimv C=\dimv I_{B_-}.
\end{equation}
For this, we
first note that we have the three triangles
\[
\xymatrix{N \ar[r]^{H^0(i)} & H^0 B_- \ar[r] & \cok(H^0 i)\oplus SK
\ar[r] & SN }
\]
\[
\xymatrix{SP_{B_-}\ar[r] & B_- \ar[r] & H^0 B_- \ar[r]^0 & I_{B_-}}
\]
\[
N \to Y \to M \to SN
\]
in $\cc_Q$. Note that $H^0 i$ is the composition of the morphism
$N \to B_-$ with the projection $B_- \to H^0 B_-$. If we form the
octahedron associated with this composition, the three triangles
we have just mentioned appear among its faces, as well as a
new triangle, namely
\[
\xymatrix{ SP_{B_-} \ar[r] & M \ar[r] & \cok(H^0 i)\oplus SK \ar[r] &I_{B_-}}.
\]
 If we
apply $H^*$ to this triangle, we obtain the exact
sequence of $H$-modules
\[
\xymatrix{0 \ar[r] & M \ar[r] & \cok(H^0 i)\oplus H^0(\tau K)
\ar[r] & I_{B_-} \ar[r] & H^0(\tau M)}.
\]
Since $M$ is an indecomposable module, $\tau M$ is either an
indecomposable non injective module or zero. The image of
$I_{B_-} \to \tau M= H^0 \tau M$ is injective (as a quotient
of an injective module). Hence it is zero and we get an exact
sequence
\[
\xymatrix{0 \ar[r] & M \ar[r] & \cok(H^0 i)\oplus H^0(\tau K)
\ar[r] & I_{B_-} \ar[r] & 0}.
\]
In the Grothendieck group, this yields
\[
0 = \dimv M -\dimv \cok(H^0 i) - \dimv H^0(\tau K) + \dimv I_{B_-} = \dimv C - \dimv H^0(\tau K) +\dimv I_{B_-}.
\]
Now, by the third triangle, $K$ is a quotient of $H^0(S^{-1}
M)=H^0(\tau^{-1}M)$. As $M$ is a non injective indecomposable module,
$H^0(\tau^{-1}M)=\tau^{-1}M$, so $\tau K$ is a quotient of $M$, and
hence, $\tau K$ is a module. Thus, we get formula~\ref{wanted} as
desired. This ends the proof of theorem~\ref{exchange}.
\end{subsection}

\section{A denominator theorem.}\label{}

\begin{subsection}{Weakly positive Laurent polynomials}
We recall an idea from \cite{BMR2}: Define an integer polynomial $P$
in $n$ variables $x_1,\ldots, x_n$ to be {\em weakly positive} if we
have $P(z)>0$ for each point $z$ of $\N^n$ which has at most one
vanishing component. If $P$ is an integer polynomial and $d$ an
$n$-tuple of integers, then $L=P/x^d$ is a {\em Laurent polynomial},
where $x^d$ is the product of the $x_i^{d_i}$.  It is a {\em weakly
positive Laurent polynomial} if we can choose $P$ weakly
positive. Clearly, in this case, no factor $x_i$ divides $P$ so that
the factorization $L=(1/x^d)\cdot P$ is unique. We call $x^d$ the {\em
denominator} of $L$. For example, the Laurent polynomial $L=x_1$ is
weakly positive with denominator $1/x_1$.  The proof of the following
lemma is elementary and left to the reader.

\begin{lemma}
\begin{itemize}
\item[a)] If $L_1$ and $L_2$ are weakly positive Laurent polynomials,
so is their sum $L_1+L_2$. Moreover, if $x^d$ and $x^e$ are the
denominators of $L_1$ and $L_2$, the denominator of $L_1+L_2$ is
$x^{\sup(d,e)}$, where by definition $\sup(d,e)_i=\sup(e_i,d_i)$,
$1\leq i\leq n$.
\item[b)] Suppose that $L_1$ and $L_2$ are Laurent polynomials and
$L_1$ is weakly positive. Then $L_2$ is weakly positive iff $L_1 L_2$
is weakly positive. Moreover, if $x^d$ and $x^e$ are the denominators
of $L_1$ and $L_2$, the denominator of $L_1 L_2$ is $x^{d+e}$.
\end{itemize}
\end{lemma}

\end{subsection}

\begin{subsection}{Denominators and dimension vectors}
From the multiplication formula, we obtain the following
denominator property for exceptional modules. A direct
proof for arbitrary modules has recently been obtained
in \cite{Hubery06}.

\begin{theorem}\label{}
Let $M$ be an indecomposable {\em exceptional} $H$-module with
dimension vector $\dimv M=(m_i)$. Then the denominator of $X_M$ as
an irreducible fraction of integral polynomials in the variables
$x_i$ is $\prod_ix_i^{m_i}$.
\end{theorem}

\begin{proof} Let us start with some preliminary remarks:
By the explicit formula for $X_M$, its denominator as an irreducible
fraction of integral polynomials is a monomial
\[
x^{\den(M)}= \prod_i x_i^{\den(M)_i} \ko
\]
where $\den(M)\in\Z^n$.  We claim that for each exceptional indecomposable
$M$, the Laurent polynomial $X_M$ is weakly positive.
Indeed, if $(T,T')$ is an exchange pair of exceptional
objects (\cf \cite{BMRRT}) and
\[
T \to B \to T' \to ST \mbox{ and } T' \to B' \to T \to ST'
\]
are non split triangles, then we have
\[
X_{T'} = \frac{X_B + X_{B'}}{X_T}
\]
by the multiplication formula. Thus, by the lemma above, if
$X_B$, $X_{B'}$ and $X_T$ are weakly positive, so is $X_{T'}$.
The claim therefore follows from the facts that each
exceptional object is a direct summand of a cluster tilting object
and that the cluster tilting graph is connected, \cf \cite{BMRRT}.
The lemma also shows that for an exchange pair $(T,T')$, we have
\[
\den(T) + \den(T') = \sup(\den(B),\den(B')).
\]
Now for an object $X$ of $\cc_Q$, we define
\[
\delta(X)=\dimv H^0(X) - \sum m_i e_i \ko
\]
where $m_i$ is the multiplicity of $SP_i$ in the decomposition of
$X$ into indecomposables and $e_i$ the $i$th vector of the
canonical basis of $\Z^n$. With this notation, we have to prove
that
\[
\delta(M)=\den(M)
\]
for each indecomposable exceptional $M$. If $P$ is indecomposable
projective, this equality is trivial for $M=SP$ and is shown for
$M=P$ by computing $X_P$ explicitly as in \cite[lemma 3.9]{caldchap}. To
prove that the equality holds for every indecomposable
exceptional, it suffices therefore, by induction, to prove that we have
\[
\delta(T) + \delta(T') = \sup(\delta(B),\delta(B')).
\]
if $(T,T')$ is an exchange pair and $B$, $B'$ are the non split
extensions of $T'$ by $T$ and $T$ by $T'$. This will be proved by
distinguishing two cases according to whether $T$ and $T'$ are
modules or one of them is a shifted projective (note that they
cannot both be shifted projectives since $\Ext^1_{\cc_Q}(T,T')\neq
0$).

{\em Case 1: $T$ and $T'$ are modules.} Since we have
\[
k=\Ext^1_{\cc_Q}(T_1,T_1')=\Ext^1_{kQ}(T_1, T_1')\oplus D\Ext^1_{kQ}(T_1',T_1),
\]
exactly one of the triangles
\[
T \to B \to T' \to ST \mbox{ and } T' \to B' \to T \to ST'
\]
comes from an exact sequence of modules. Let us assume it is the
first one. Then, by applying the functor $H^0$, we get exact
sequences
\[
0 \to T \to B \to T' \to 0 \mbox{ and } T' \to H^0(B') \to T.
\]
These show that we have
\[
\dimv B = \dimv T + \dimv T' \mbox{ and } \dimv H^0(B')\leq \dimv T
+\dimv T' \leq \dimv B.
\]
It follows that $\delta(B')\leq \dimv B = \delta(B)$ and
\[
\delta(T')+\delta(T)=\delta(B) = \sup(\delta(B), \delta(B')).
\]

{\em Case 2: We have $T=SP$ for an indecomposable projective $P$
and $T'$ is a module.} Again, we have non split triangles
\[
T \to B \to T' \to ST \mbox{ and } T' \to B' \to T \to ST.
\]
We have
\[
k=\Ext^1_{\cc_Q}(T, T')= \Hom_{\cc_Q}(SP, ST')= \Hom_{kQ}(P, T')
\]
and
\[
k=\Ext^1_{\cc_Q}(T', T)= \Hom_{\cc_Q}(T', S^2 P)=
\Hom_{kQ}(T',\nu P) \ko
\]
where $\nu$ is the Nakayama functor. Since $\mod kQ$ is hereditary,
if $f: L\to M$ is a morphism of $\mod kQ$, then in the derived
category, we have a triangle
\[
\xymatrix{ L \ar[r]^f & M \ar[r] & \cok(f)\oplus S\ker(f) \ar[r] & SL.}
\]
It follows that the triangles above are in fact isomorphic to the
triangles
\[
S^{-1}\nu P \to S^{-1}\cok(f)\oplus \ker(f) \to T' \to \nu P
\mbox{ and }
T' \to S\ker(g)\oplus \cok(g) \to SP \to ST'
\]
associated with arbitrary non zero morphisms
$f:T' \to \nu P$ and $g: P \to T'$.
Now $\cok(f)$ is injective as a quotient of an injective
and $\ker(g)$ projective as a quotient of a projective.
Moreover, if $i$ is the vertex of $Q$ such that $P=P_i$, then,
as a submodule of $P_i$, the module $\ker(f)$ is a direct sum of
indecomposables $P_j$ such that $j<i$. Similarly, as a
quotient of $\nu P_i$, the module $\cok(g)$ is a direct sum
of indecomposables $\nu P_j$ such that
$i<j$ (note that we consider right modules and order the vertices of
$Q$ in the natural way). It follows that we have
\[
\sup(\delta(S^{-1}\cok(f)), \delta(S\ker(g)))=0
\]
(note that both vectors have negative components).
Thus, we have
\begin{align*}
\sup(\delta(B),\delta(B')) & = \sup(\dimv \ker(f) + \delta(S^{-1} \cok(f)), 
\dimv \cok(g) + \delta(S\ker(g))) \\
& =\sup(\dimv \ker(f), \dimv \cok(g)).
\end{align*}
It remains to be proved that
\begin{equation}\label{dimvt}
\dimv T' + \delta(T)= \sup(\dimv \ker(f), \dimv \cok(g)). \quad\quad
\end{equation}
We check this equality by comparing both sides at each vertex $j$
of $Q$. As above, let $i$ be the vertex of $Q$ such that $P=P_i$ so that
we have $\delta(T)=-e_i$. We have
\[
\Hom(P_i, T') = \Ext^1_{\cc_Q}(SP_i, T') = k
\]
so that $(\dimv T')_i=1$. The maps $f: T' \to \nu P$ and $g:P\to
T'$ induce isomorphisms in $\Hom(P_i, ?)$ since $g\circ f$
induces an isomorphism between the one-dimensional
spaces $\Hom(P_i, P)$ to $\Hom(P_i, \nu P)$. It follows that
$(\dimv \ker(f))_i$ and $(\dimv \cok(f))_i$ both vanish so that the
equality (\ref{dimvt}) holds at $j=i$. Now consider the exact sequences
\[
\xymatrix{ 0 \ar[r] & \ker(f) \ar[r] & T' \ar[r]^{f} & \nu P}
\mbox{ and }
\xymatrix{ P \ar[r]^{g} & T' \ar[r] & \cok(g) \ar[r] & 0.}
\]
Suppose that $j$ is not a
predecessor of $i$. Then $(\dimv P)_j=0$ and we have
$(\dimv T')_j = (\dimv \cok(g))_j$ by the second sequence and
$(\dimv T')_j \geq (\dimv \ker(f))_j$ by the first so that
\[
(\dimv T')_j = (\dimv \cok(g))_j = \sup((\dimv \cok(g))_j, (\dimv
(ker(g)))_j)
\]
and (\ref{dimvt}) holds at $j$. Similarly, if $j$ is not a successor of
$i$, we see that the equality (\ref{dimvt}) holds at $j$. Since $Q$ has no
oriented cycles, each vertex $j\neq i$ of $Q$ is a non-successor
or a non-predecessor of $i$. Thus, the proof of (\ref{dimvt}) is complete.
\end{proof}
\end{subsection}

\begin{subsection}{}
As a corollary of the denominator theorem, we will prove an
injectivity property of the map $M\mapsto X_M$.

We recall first a few facts on quiver representations.

A representation of $Q$ over a field $F$ is a $Q_0$-graded
$F$-vector space $V=\oplus_{i\in Q_0} V_i$ together with an element
$x=(x_h)_{h\in Q_1}$ in $E_V:=\prod_{h\in Q_1}\Hom(V_{s(h)},V_{t(h)})$,
where $s(h)$ is the source and $t(h)$ the target of
the arrow $h$. The group $G_V:=\prod_{i\in Q_0} \GL(V_i)$ acts on $E_V$ by
$(g_i).(x_h)=(g_{t(h)}x_hg_{s(h)}^{-1})$. A representation $(M,x)$
over a field $F$ can be functorialy considered as an $FQ$-module and the dimension
vector of this module is $\dimv M=(\dim M_i)$.

Clearly, the isoclasses of
finite-dimensional $FQ$-modules are naturally identified with
$G_V$-orbits of representations of $Q$.

\begin{corollary}\label{injective}
If $M$ and $M'$ are non isomorphic indecomposable modules without
self-extensions, then $X_M\not=X_{M'}$.
\end{corollary}

\begin{proof}
It is well known that in the identification above, an isoclass of
$kQ$-module with no self-extension corresponds to an orbit which is
dense in its representation space $E_V$.
Therefore, if $M$ and $M'$ are non isomorphic modules without self-extensions,
the corresponding orbits cannot be in the same representation
space. Hence, $M$ and $M'$ cannot have the same dimension vector.

By the theorem above, we conclude that $X_M\not=X_{M'}$.

\end{proof}

\end{subsection} 

\section{Application to a class of cluster algebras}\label{}
\begin{subsection}{}
We recall some terminology on cluster algebras.  The reader can find
more precise and complete information in \cite{fominzelevinsky2}.\par

Let $n$ be a positive integer. We fix the {\it ambient field} ${\cF}=
\Q(x_1,\ldots,x_n)$, where the $x_i$'s are indeterminates. Let
${\mathbf x}$ be a free generating set of $\cF$ over $\Q$ and let
$B=(b_{ij})$ be an $n\times n$ antisymmetric matrix with coefficients
in $\Z$. Such a pair $({\mathbf x},B)$ is called {\it a seed}.\par

Let $({\mathbf u}, B)$ be a seed and let $u_j$, $1\leq j\leq n$,
be in ${\mathbf u}$.  We define a new seed as follows. Let $u_j'$ be
the element of $\cF$ defined by the {\it exchange relation}:
\begin{equation}\label{exchangerelation}
u_ju_j'=\prod_{b_{ij}>0}u_i^{b_{ij}}+\prod_{b_{ij}<0}u_i^{-b_{ij}}.
\end{equation}
Set ${\mathbf u'}={\mathbf u}\cup\{u_j'\}\backslash \{u_j\}$.  Let
$B'$ be the $n\times n$ matrix given by
\[
b_{ik}'=\begin{cases}-b_{ik}&\hbox{if } i=j \hbox{ or } k=j\\
b_{ik}+\frac{1}{2}( \,|b_{ij}|\, b_{jk}+b_{ij}\, |b_{jk}|\,) & \hbox{
otherwise.}\cr\end{cases}
\]
By a result of Fomin and Zelevinsky, $({\mathbf u'},B')=\mu_j({\mathbf
u}, B)$ is a seed.  It is called the {\it mutation} of the seed
$({\mathbf u},B)$ in the direction $u_j$ (or $j$).  We consider all
the seeds obtained by iterated mutations. The free generating sets
occurring in the seeds are called {\it clusters}, and the variables
they contain are called {\it cluster variables}. By definition, the
{\it cluster algebra} $\ca({\mathbf x}, B)$ associated to the seed
$({\mathbf x},B)$ is the $\Z$-subalgebra of $\cF$ generated by the set
of cluster variables. The graph whose vertices are the seeds and whose
edges are the mutations between two seeds is called the {\it mutation
graph} of the cluster algebra.

The {\it Laurent phenomenon}, see \cite{fominzelevinsky1},
asserts that the cluster variables are Laurent polynomials with integer
coefficients in the $x_i$, $1\leq i\leq n$.  So, we have $\ca({\mathbf
x},B)\subset\Z[x_1^{\pm 1},\ldots,x_n^{\pm 1}]$.

Note that an antisymmetric matrix $B$ defines a quiver $Q=Q_B$ with
vertices corresponding to its rows (or columns) and which has $b_{ij}$
arrows from the vertex $i$ to the vertex $j$ whenever $b_{ij}\geq 0$.
The quivers $Q$ thus obtained are precisely the finite quivers
without oriented cycles of length $1$ or $2$. For such quivers $Q$, we
denote by $B_Q$ the corresponding antisymmetric matrix. The cluster
algebra associated to the seed $({\mathbf x},B)$ will be also denoted
by $\ca(Q)$. In the sequel, we will be concerned with cluster algebras
associated to a quiver $Q$ without oriented cycles.
\end{subsection}
\begin{subsection}{}
We fix a quiver $Q$ without oriented cycles and we set $H=kQ$. We
consider the cluster category $\cc=\cc_H$ associated to the quiver
$Q$, \cf \cite{BMRRT}. An object $T$ of $\cc$ is called {\it
exceptional} if it has no self-extensions, \ie if $\Ext^1(T,T)=0$.
An exceptional object is called {\it cluster-tilting} or simply
{\it tilting} (although this is an abuse of language) if it has
$n$ non isomorphic indecomposable direct summands, where $n$ is
the number of vertices of $Q$. In the sequel, we will often
identify a tilting object with the datum of its indecomposable
summands. An exceptional object is called {\it almost tilting} if
it has $n-1$ non isomorphic indecomposable direct summands. It was
shown in \cite{BMRRT} that any almost tilting object $\overline T$
can be completed to precisely two non isomorphic tilting objects
$T$ and $T^*$.

For any tilting object $T$ of $\cc$, let $Q_T$ be the quiver
associated to the algebra $\End_{\cc}(T)$.  To be explicit, fix an
ordering of the indecomposable summands $T_1,\ldots,T_n$ of $T$ and let
$A$ be the endomorphism algebra of the sum of the $T_i$.  Let $e_i \in
A$ the idempotent corresponding to $T_i$. Then the vertices of $Q_T$
are $1,\ldots,n$, and the number of arrows from $i$ to $j$ is equal to
$\dim e_j((\rad A)/(\rad A)^2) e_i$.  A pair $(T,Q_T)$ is called
a {\it cluster seed}.

For $1\leq i\leq n$, we define, following \cite{BMR2}, the mutation of the
cluster seed $(T,Q_T)$ in direction $i$ by
\[\delta_i(T,Q_T):=(T^*,Q_{T^*}),\]
where $T$ and $T^*$ are the two completions of the almost tilting object
\[\overline T=T_1\oplus\ldots \oplus T_{i-1}\oplus T_{i+1}\oplus\ldots
\oplus T_n.\]
Note that there exists an indecomposable object $T_i^*$, unique
up to isomorphism, such that
\[T^*=T_1\oplus\ldots \oplus T_{i-1}\oplus T_i^*\oplus T_{i+1}\oplus\ldots
\oplus T_n,\]
which provides a natural ordering of the indecomposable summands of $T^*$.

The following theorem is the main result of this article. The first
assertion is a refinement of Conjecture~9.1 of \cite{BMRRT} and the
second assertion strengthens the main result of \cite {BMR2}.
\begin{theorem}\label{main}
Let $Q$ be a quiver with $n$ vertices and no oriented cycles, and let
$H=kQ$ be the hereditary algebra associated to $Q$. Then
\begin{itemize}
  \item[(i)] The correspondence $M\mapsto X_M$ provides a bijection
  between the set of indecomposable objects without self-extensions
  of $\cc_H$ and the set of cluster variables of $\ca(Q)$.
  \item[(ii)] The correspondence
  $\{T_1,\ldots,T_n\}\mapsto\{X_{T_1},\ldots,X_{T_n}\}$ provides a
  bijection compatible with mutations between the set of tilting
  objects of $\cc_H$ and the set of clusters of $\ca(Q)$.
\end{itemize}
\end{theorem}
\begin{proof}
By construction, any cluster variable belongs to a cluster. As the
map $M\mapsto X_M$ is injective on the set of indecomposable
objects without self-extensions by corollary~\ref{injective}, it
is enough to prove (ii).

Let us prove (ii). Suppose that $T=T_1\oplus\ldots \oplus T_n$ is a
tilting object of $\cc$ and let $T^*$ be its mutation in direction
$i$. Then $\Ext^1(T_i,T_i^*)$ is one-dimensional
by \cite{BMRRT}. Hence, by theorem \ref{exchange}, we have
\begin{equation}\label{tiltingexchange}
X_{T_i}X_{T_i^*}=\prod_jX_{T_j}^{a_{ij}}+\prod_jX_{T_j}^{c_{ij}},
\end{equation}
where $a_{ij}$ and $c_{ij}$ are integers defined by the following non
split triangles (unique up to isomorphism)
\[T_i\rightarrow\oplus a_{ij}T_j\rightarrow T_i^*\rightarrow ST_i\]
\[T_i^*\rightarrow\oplus c_{ij}T_j\rightarrow T_i\rightarrow ST_i^*.\]
By a theorem~6.2 b) of \cite{BMR2}, the quiver $Q_T$
is determined by these triangles: for any $i$ and $j$, there are
$a_{ij}$ arrows from $i$ to $j$ and $c_{ij}$ arrows from $j$ to
$i$. Moreover, if there exists an arrow from $i$ to $j$ , then
there is no arrow from $j$ to $i$, by proposition~3.2 of \cite{BMR2}.

We now define, as in \cite{BMR2}, a correspondence $\beta$ between
tilting seeds and cluster seeds.  First note that the shift of $H$,
is a tilting object and that $(SH,Q)$ is a tilting seed.
For a given word $i_1\ldots i_t$, we can define
\begin{equation}\label{beta0}
\beta(SH,Q)=(\mathbf x, B_Q),
\end{equation}
\begin{equation}\label{beta1}
\beta(\delta_{i_t}\ldots\delta_{i_1}(SH,Q))=\mu_{i_t}\ldots\mu_{i_1}(\mathbf x,B_Q).
\end{equation}
Set $(T,Q_T):=\delta_{i_t}\ldots\delta_{i_1}(SH,Q))$.  By \cite{BMR2},
the quiver obtained from $Q$ by the sequence of {\it tilting} mutations
in the directions $i_1,\ldots,i_t$ is equal to the quiver obtained from $Q$
by the sequence of {\it cluster} mutations in the directions
$i_1,\ldots,i_t$. Hence, by comparing the cluster exchange relation
\ref{exchangerelation} and the tilting exchange relation
\ref{tiltingexchange}, we obtain by induction that
\[ \beta(\delta_{i_t}\ldots\delta_{i_1}(SH,Q))=
(\{X_{T_1},\ldots,X_{T_n}\}, B_{Q_T}).\]
In particular, $\beta(\delta_{i_t}\ldots\delta_{i_1}(SH,Q))$ does not
depend on the choice of the word $i_1\ldots i_t$.

By proposition~3.5 of \cite{BMRRT}, the mutation graph on the set of
tilting seeds is connected. Hence, equalities~\ref{beta0} and
\ref{beta1} define a map $\beta$ from the complete set of
tilting seeds to the set of cluster seeds. The surjectivity of $\beta$
follows from the fact that its image is stable under mutation. The
injectivity of $\beta$ follows from corollary \ref{injective}.
\end{proof}

\end{subsection}
\begin{subsection}{} \label{ss:connectedness}
This section is devoted to the proof of some of the conjectures
formulated by S.~Fomin and A.~Zelevinsky in \cite{FZ-conf}.  
The first corollary is a straightforward consequence of theorem~\ref{main}.
It corresponds to \cite[Conjecture 4.14 (2)]{FZ-conf} in the acyclic case.
\begin{corollary}
Let $Q$ be an finite quiver without oriented cycles. Then a
cluster seed $(\mathbf u,B)$ of $\ac(Q)$ only depends on $\mathbf u$.
\end{corollary}
This corollary is  \cite[Conjecture 4.14 (3)]{FZ-conf} in the acyclic case.

\begin{corollary}
For any cluster variable $x$, the set of seeds whose clusters contain $x$
form a connected subgraph of the exchange graph.
\end{corollary}
\begin{proof} Indeed, the cluster variable $x$ corresponds to an
exceptional indecomposable object $T_1$ of $\cc_Q$. Without
restriction of generality, we assume that $T_1$ is non projective.
The seeds containing $x$ are in bijection with the {\em
completions of $T_1$}, \ie the sets $\{T_2, \ldots, T_n\}$ of
indecomposables such that the sum of the $T_i$ is cluster tilting.
Two seeds are joined by an edge of the exchange graph
iff the corresponding sets of exceptional indecomposables are
obtained from one each other by a mutation. By
\cite{BMRRT}, this occurs iff they differ by precisely two
indecomposables $T_i$ and $T_i^*$ and these satisfy
\[
\dim\Ext^1(T_i, T_i^*)=1.
\]
This makes it clear that theorem~\ref{thm:clusterreduction} below
yields a bijection compatible with mutations
\[
\{T_2, \ldots, T_n\} \mapsto \{PT_2, \ldots, PT_n\}
\]
between the completions of $T_1$ and the basic tilting sets of
$\cc_{Q'}$, where $Q'$ is the quiver of the endomorphism ring of a
projective generator of the category $\ch'\subset \mod kQ$ of
modules $L$ with
\[
\Hom(M,L)=0=\Ext^1(M,L).
\]
Thus, by theorem~\ref{main} (ii),
the subgraph of the exchange graph of $Q$ formed by the seeds
containing $x$ is isomorphic to the exchange graph of $Q'$, which
is connected by definition.
\end{proof}
A consequence of theorem \ref{main} is also the proof of
 \cite[Conjecture 4.14 (4)]{FZ-conf} in the general case.

\begin{corollary}
The set of seeds whose matrix is acyclic form a connected subgraph
(possibly empty) of the exchange graph.
\end{corollary}
\begin{proof} A seed with an acyclic matrix corresponds to a
cluster tilting object $T$ whose endomorphism algebra
$A=\End_{\cc_Q}(T)$ has a quiver without oriented cycles. 
Thus, the algebra $A$ is both, cluster-tilted and of
finite global dimension. By corollary~2.1 of \cite{KellerReiten05},
it is hereditary. So the category $\cc_A$ is well-defined
and the equivalence between the derived categories of
$A$ and $Q$ induces a triangle equivalence $\cc_A \iso \cc_Q$ which takes
$A$ to $T$. Such an equivalence induces an isomorphism
\[
\Gamma_A \to \Gamma_B
\]
of the Auslander-Reiten quivers of the two cluster categories. We
refer to \cite{BMRRT} for the description of the Auslander-Reiten
quivers.  Since $A$ is hereditary, the quiver of its
indecomposable projectives forms a slice of the component
$\Gamma_A^{pr}$ of $\Gamma_A$ containing the projectives (recall
that a slice is a full connected subquiver whose vertices are a
system of representatives of the $\tau$-orbits in the component).
The isomorphism must take $\Gamma_A^{pr}$ to $\Gamma_B^{pr}$ since
this is the only components isomorphic to the repetition $\Z R$ of
a finite quiver $R$. It is clear that any slice of $\Gamma_B^{pr}$
can be transformed to the slice of the projectives by finitely
many reflections at sources or sinks.
\end{proof}

\end{subsection}

\subsection{Cluster tilting objects containing a given summand}
\label{ss:fixed-summand} Here, we refine a technique pioneered in
section~2 of \cite{BMR2}: Let $H$ be a finite-dimen\-sional
hereditary algebra and $\ch$ the category of finite-dimen\-sional
right $H$-modules. Let $M \in\ch $ be a non projective
indecomposable with $\Ext^1(M,M)=0$. Then $\End(M)$ is a (possibly
non commutative) field. Let $\ch'$ be the full subcategory on the
modules $L$ such that
\[
\Hom(M,L)=0 \mbox{ and } \Ext^1(M,L)=0.
\]
We know from \cite{HappelRickardSchofield88} and \cite{happel}
that $\ch'$ is a hereditary abelian category with enough projectives and
that a projective generator $G$ of $\ch'$ is obtained by choosing
an exact sequence
\[
0 \to H \to G \to M^r \to 0
\]
which induces an isomorphism
\[
\Hom(M,M^r) \iso \Ext^1(M,H).
\]
Let $\cc_\ch$ and $\cc_{\ch'}$ be the cluster categories
associated with $\ch$ and $\ch'$. The following theorem
is an elaboration on theorem 2.13 of \cite{BMR2}.

\begin{theorem} \label{thm:clusterreduction}
Let $\cc(\ch,M)$ be the full additive subcategory of $\cc_\ch$
whose objects are the sums of indecomposables $L$ of $\cc_\ch$
such that $\Ext^1(M,L)=0$. There is a canonical equivalence
of $k$-linear categories
\[
P:\xymatrix{\cc(\ch,M)/(M) \ar[r]^-{\sim} & \cc_{\ch'} }\ko
\]
where $(M)$ denotes the ideal of morphisms factoring through
a sum of copies of $M$. Moreover, we have
\[
\Ext^1(L_1, L_2) \cong \Ext^1(PL_1, PL_2)
\]
for all $L_1, L_2 \in \cc(\ch,M)$.
\end{theorem}

Note that $\cc(\ch,M)$ is not a triangulated subcategory and not
even stable under the shift functor. The theorem merely claims
that as a $k$-linear category, $\cc_{\ch'}$ is a `subquotient' of
$\cc_\ch$. To construct the equivalence $P$, we choose a
`fundamental domain' for the action of the autoequivalence
$F=\tau^{-1}S$ on $\cd$.

Let $\cp$ be the full subcategory of the projectives of $\ch$ and
$\ch^+$ the full additive subcategory of $\cd=\cd^b(\ch)$ each of
whose indecomposables lies in $\ch$ or $S\cp$. Let $\pi: \cd \to
\cc_\ch$ be the projection functor. We know from \cite{BMRRT} that
$\pi$ induces a bijection from the set of isoclasses of
indecomposables of $\ch^+$ to that of $\cc_\ch$ and that we have
\[
\Ext^1(\pi(L_1),\pi(L_2)) \iso \Ext^1(L_1, L_2) \oplus
D\Ext^1(L_2, L_1)
\]
for any two indecomposables of $\ch^+$. Moreover, the category
$\cc_\ch$ is equivalent to the category whose objects are those
of $\ch^+$ and whose morphisms are given by
\[
\Hom(L_1, L_2) \oplus \Hom(L_1, F L_2)
\]
with the natural composition. Therefore,
theorem~\ref{thm:clusterreduction} follows from

\begin{theorem} There is a canonical bijection $L\mapsto L'$ from
the set of isoclasses of indecomposables $L$ of $\ch^+$ with
\[
L\not\cong M \ko \Ext^1(L,M)=0 \mbox{ and } \Ext^1(M,L)=0 \quad\quad\quad (*)
\]
to the set of isoclasses of indecomposables of $\ch'^+$. Moreover, for any two
objects $L_1, L_2$ of $\ch^+$ satisfying $(*)$, there is a
canonical isomorphism
\[
\Ext^1(L_1,L_2) \iso \Ext^1(L_1', L_2')
\]
and there are canonical isomorphisms
\[
\Hom(L_1,L_2)/(M) \iso \Hom(L_1', L_2') \mbox{ and }
\Hom(L_1, F L_2)/(M) \iso \Hom(L_1', F'  L'_2)
\]
compatible with compositions.
\end{theorem}

Before giving the proof, let us illustrate the statement
on the following example: We consider the path algebra $H=kQ$ of
a linearly oriented quiver $Q$ of type $A_6$. Below, we have
drawn the Auslander-Reiten quiver of its derived category.
The vertices corresponding to indecomposables concentrated in
degree $0$ lie between the two hatched lines.
We use the symbols
\begin{itemize}
\item[$\sq$] for the $14$ indecomposables $L$ of $\ch^+$ not isomorphic to
$M$ and which satisfy
\[
\Ext^1(L,M)=0=\Ext^1(M,L)\ko
\]
\item[$\bt$] for the $5$ indecomposable projectives of $\ch'$,
\item[$\point$] for indecomposable non projective objects of
$\ch'$,
\item[$\dmd$] for shifted copies $SP$ of projectives of $\ch'$.
\end{itemize}
Notice the two rectangular zones starting from $S^{-1}M$
respectively ending in $SM$ where no $\sq$ occurs.
If $L\mapsto L'$ denotes the map of the theorem, we have
a triangle
\[
L_\cu \to L \to L' \to SL_\cu
\]
where $L_\cu$ is a sum of copies of $M$ and $L_\cu \to L$ induces
a bijection $\Hom(M,L_\cu) \iso \Hom(M,L)$. Thus we have $L=L'$
for all $L$ with $\Hom(M,L)=0$. The corresponding triangles
for the others are visible in the diagram below.
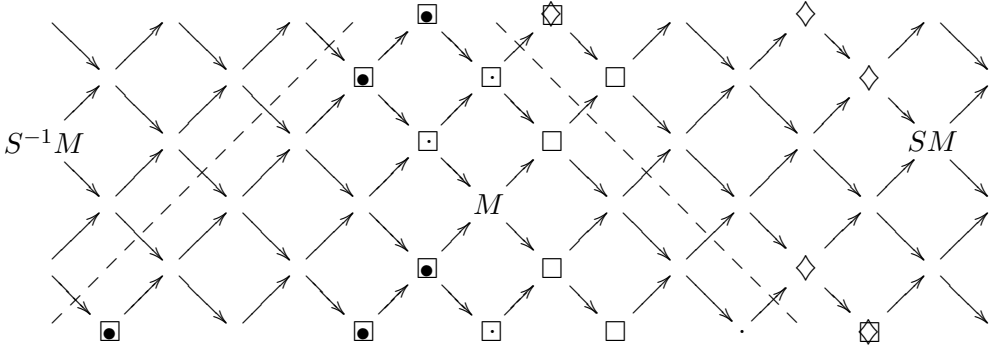
\[
\xymatrix@!0{
\ar[rd] & & \ar[rd] & & \ar[rd] & & \sqbt \ar[rd] & \ar@{--}[rrrrrddddd]
& \sqdmd \ar[rd] & &
\ar[rd] & & \dmd \ar[rd] & & \ar[rd] \\
 & \ar[ru] \ar[rd] & & \ar[ru] \ar[rd] & &
\sqbt\ar[ru] \ar[rd] & & \sqpt \ar[ru] \ar[rd] & & \sq \ar[ru] \ar[rd] & &
\ar[ru] \ar[rd] & & \dmd \ar[ru] \ar[rd] & & \\
S^{-1}M \ar[ru] \ar[rd] & & \ar[ru] \ar[rd] & & \ar[ru] \ar[rd] & &
 \sqpt  \ar[ru] \ar[rd] & & \sq \ar[ru] \ar[rd] & & \ar[ru] \ar[rd] & &
\ar[ru]  \ar[rd] & & SM \ar[ru] \ar[rd] \\
& \ar[ru]\ar[rd] & & \ar[ru] \ar[rd]  & & \ar[ru] \ar[rd] & &
M \ar[ru] \ar[rd] & & \ar[ru] \ar[rd] & & \ar[ru] \ar[rd] & &
\ar[ru] \ar[rd] & & \\
\ar[ru] \ar[rd] & & \ar[ru] \ar[rd] & & \ar[ru] \ar[rd] & &
\sqbt \ar[ru] \ar[rd] & & \sq \ar[ru] \ar[rd] & & \ar[ru] \ar[rd] & &
\dmd \ar[ru] \ar[rd] & & \ar[ru] \ar[rd] & \\
\ar@{--}[rrrrruuuuu] & \sqbt \ar[ru] & & \ar[ru] & &
\sqbt \ar[ru] & & \sqpt \ar[ru] & & \sq \ar[ru] & &
\point \ar[ru] & & \sqdmd\ar[ru] & &
}
\]

Several of the arguments needed in the proof are contained
in section~2 of \cite{BMR2}. For the convenience of the
reader, we nevertheless include them below.

\begin{proof} Let $\cu\subset \cd$ be the full triangulated
subcategory generated by $M$. Since $\Ext^1(M,M)$ vanishes and
$\Hom(M,M)$ is a field, its objects are the sums of shifted copies
of $M$. Let $\cv$ be the full subcategory of $\cd$ whose objects
are the $L\in \cd$ such that $\Hom(U, L)=0$ for all $U\in\cu$.
Then $\cu,\cv$ form a semiorthogonal decomposition
\cite{BondalKapranov89} of $\cd$, \ie
for each object $X$ of $\cd$, there is a triangle
\[
X_\cu \to X \to X^\cv \to SX_\cu
\]
with $X_\cu\in \cu$ and $X^\cv\in \cv$. This triangle is unique up to
unique isomorphism; the functor $X\to X_\cu$ is right adjoint to the
inclusion of $\cu$ and the functor $X \mapsto X^\cv$ is left adjoint
to the inclusion of $\cv$. We have $\ch'=\ch\cap \cv$ and the
inclusion $\ch' \subset \cv$ extends canonically to an equivalence
$\cd^b(\ch') \to \cv$. In particular, each object of $\cv$ is a direct
sum of shifts of objects of $\ch'$. We have $\cu\cap\ch = \cm$, the
full subcategory on the direct sums of copies of $M$. The inclusion
$\ch'\subset \ch$ commutes with kernels, cokernels and preserves
$\Ext^1$-groups. We will show that $L \mapsto L'=L^\cv$ yields the
bijection announced in the assertion.

Let $L$ be indecomposable in $\ch^+$ such that $(*)$ holds.
Let us first show that $\Hom(S^i M,L)$ vanishes if $i\neq 0$.
Indeed, if $L$ belongs to $\ch$, then this group clearly vanishes
if $i\neq 0,-1$ and if $i=-1$, it vanishes because
$\Ext^1(M,L)=0$. If $L=SP$ for a projective $P\in\ch$, then
$\Hom(S^iM,L)=\Hom(S^i M, SP)$ clearly vanishes for $i\neq 0,1$ and
it vanishes for $i=1$ because $M$ is a non projective indecomposable.

Now let us show that $L^\cv$ is indecomposable:
Consider the canonical triangle
\[
L_\cu \to L \to L^\cv \to SL_\cu.
\]
Since $\Hom(S^i M, L)$ vanishes for $i\neq 0$, we have $L_\cu \in \cm$.
Therefore, in the associated exact sequence
\[
\Hom(L,L) \to \Hom(L,L^\cv) \to \Hom(L,SL_\cu)
\]
the third term vanishes.  Thus the composition
\[
\Hom(L,L) \to \Hom(L,L^\cv) \iso \Hom(L^\cv,L^\cv)
\]
is surjective and $\End(L^\cv)$ is local as a quotient of the
local ring $\End(L)$.

Let us show that $L'$ belongs to $\ch'^+$. Since $L_\cu$ belongs
to $\cm$, the canonical morphism $f: L_\cu \to L$ is a morphism
of $\ch$ and therefore its cone $L^\cv$ in $\cd$ is isomorphic
to $\cok(f)\oplus S\ker(f)$. Since $L^\cv$ is indecomposable,
one of the two summands vanishes. If $\ker(f)$ vanishes, then
$L^\cv$ belongs to $\ch'\subset \ch'^+$. If $\cok(f)$ vanishes,
we have to show that $\ker(f)$ is projective in $\ch'$. Now
indeed the short exact sequence
\[
0 \to \ker(f) \to L_\cu \to L \to 0
\]
induces a surjection
\[
\Ext^1_\ch(L_\cu, U) \to \Ext^1_\ch(\ker(f), U) \to 0
\]
for each $U\in \ch$. The left hand term vanishes since $L_\cu$
is a sum of copies of $M$ and the right hand term is isomorphic
to $\Ext^1_{\ch'}(\ker(f),U)$ because the inclusion $\ch'\subset\ch$
preserves extension groups. Thus, $\ker(f)$ is projective in $\ch'$.

From what we have shown, we conclude that the map $L \to L'$ is
well-defined. Let us show that it is injective. For this, we show
that the morphism $L^\cv \to SL_\cu$ occurring in the canonical
triangle is a minimal left $S\cm$-approximation. Then $L$ is
determined up to isomorphism as the shifted cone over this
morphism. To show that $L^\cv\to SL_\cu$ is a minimal left
approximation, consider the canonical triangle
\[
L_\cu \to L \to L^\cv \to SL_\cu
\]
and the induced sequence
\[
\Hom(L,SM) \la \Hom(L^\cv, SM) \la
\Hom(SL_\cu, SM).
\]
Since $\Hom(L,SM)$ vanishes by assumption, we do get a surjection
$\Hom(L^\cv, SM) \la \Hom(SL_\cu,SM)$. If it is not minimal,
then there is a retraction $r: SL_\cu \to SM$
whose composition
with $L^\cv \to SL_\cu$ vanishes. Then $r$ extends to a retraction
$\tilde{r}: SL \to SM$. This is impossible since $L$ is indecomposable
and not isomorphic to $M$.

Let us show now that $L \mapsto L'$ is surjective. Let $N$ be indecomposable
in $\ch'^+$. Let $N \to SM'$ be a minimal $S\cm$-approximation and
form the triangle
\[
M' \to L \to N \to SM'.
\]
Let us show that $L$ is indecomposable.
Since $M'\in\cu$, we have $L^\cv\iso N^\cv$ and since $N\in \cv$,
we have $L_\cu\iso M'$. If $L$ is decomposable,
say $L=L_1\oplus L_2$, then we get
\[
L_1^\cv \oplus L_2^\cv \iso N
\]
and, say, $L_1^\cv$ vanishes. Then $L_1$ belongs to $\cu$ and
thus $M'\iso L_1\oplus (L_2)_\cu$. Since $N \to SM'$ is a minimal
$S\cm$-approximation, we have $L_1=0$. So $L$ is indecomposable.

Let us show that $L$ belongs to $\ch^+$.
It is clear from the above triangle that $L$ has homology at most
in degrees $0$ and $1$. Since $L$ is indecomposable, its homology
is concentrated in one degree. If the homology is concentrated
in degree $0$, then $L$ belongs to $\ch\subset \ch^+$.
Suppose that $L$ has its homology concentrated in degree $1$.
Then we must have $N=SQ$ for some indecomposable
projective $Q$ of $\ch'$. We know that if $P_H$ is a projective generator for $\ch$, then
$P_H^\cv$ is is a projective generator for $\ch'$. Thus, there is
a projective $P$ of $\ch$ and a section $s: Q \to P^\cu$ which identifies $Q$ with a direct
factor of $P^\cu$. Since $N \to SM'$ is
an $S\cm$-approximation, the composition
\[
\xymatrix{ N \ar[r]^{Ss} & SP^\cu \ar[r] & SP_\cm}
\]
extends to $SM'$ so that we obtain a morphism of triangles
\[
\xymatrix{ M' \ar[r] \ar[d] & L \ar[r] \ar[d] & N \ar[r]\ar[d]^{Ss} &
SM' \ar[d] \\
SP_\cu  \ar[r] & SP \ar[r] & SP^\cu \ar[r] & SP_\cm .}
\]
The morphism $L\to SP$ is non zero since its composition with
$SP\to SP^\cu$ equals the composition of the non zero morphism
$L \to N$ with the section $Ss$. So we obtain a non zero
morphism $S^{-1}L \to P$ in $\ch$. Since $S^{-1}L$ is indecomposable
and $P$ is projective, $S^{-1}L$ is projective and we have
$L\in\ch^+$.

Finally, let us show that $L$ satisfies the condition $(*)$.
If $L$ was isomorphic to $M$, we would have $N= L^\cv =0$ contrary to our hypothesis that $N$ is indecomposable.

The triangle
\[
M' \to L \to N \to SM'
\]
yields an exact sequence
\[
\Hom(M',SM) \la \Hom(L,SM) \la \Hom(N, SM) \la \Hom(SM', SM).
\]
Here the leftmost term vanishes since $\Ext^1(M,M)=0$ and
the rightmost map is surjective since $N \to SM'$ is a left
$S\cm$-approximation. Thus we have $\Ext^1(L,M)=0$. The triangle
also yields the sequence
\[
\Hom(S^{-1}M, M') \to \Hom(S^{-1}M,L) \to \Hom(S^{-1}M, N).
\]
The left hand term vanishes since $\Ext^1(L_1,M)=0$ and
the right hand term vanishes since $N$ belongs to $\cv$.
Thus we have $\Ext^1(M,L)=0$.

Now let $L_1$, $L_2$ be indecomposables of $\ch^+$ satisfying
condition $(*)$. Consider the triangle
\[
(L_2)_\cu \to L_2 \to L_2^\cv \to S(L_2)_\cu.
\]
It induces an exact sequence
\[
\Hom(S^{-1}L_1, (L_2)_\cu) \to \Hom(S^{-1}L_1, L_2) \to
\Hom(S^{-1}L_1, L_2^\cv) \to \Hom(S^{-1}L_1, S(L_2)_\cu).
\]
The leftmost term vanishes since $\Ext^1(L_1, M)=0$ and
the rightmost term vanishes since $\Ext^2(L_1,M)=0$. Thus we
have
\[
\Hom(L_1, SL_2) \iso \Hom(L_1, SL_2^\cv) \iso \Hom (L_1^\cv, SL_2^\cv) \ko
\]
which proves the assertion on the extension groups.
The above triangle also induces an exact sequence
\[
\Hom(L_1, (L_2)_\cu) \to \Hom(L_1,L_2) \to \Hom(L_1, L_2^\cv) \to
\Hom(L_1, S(L_2)_\cu).
\]
The last term vanishes since $\Ext^1(M,M)=0$. Thus the kernel
of the map
\[
\Hom(L_1, L_2) \to \Hom(L_1, L_2^\cv) \iso \Hom(L_1^\cv, L_2^\cv)
\]
is formed by the morphisms factoring through sums of $M$.
Put $F=\tau^{-1}S$. Consider the triangle
\[
(FL_2)_\cu \to FL_2 \to (FL_2)^\cv \to S(FL_2)_\cu.
\]
Note that the functor $F$ does not take $\cv$ to itself.
We have
\[
\Hom(S^i M, FL_2) \iso D\Hom(L_2, S^i \tau^2 M).
\]
This can be non zero only if $i$ equals $0$ or $1$. Thus
$(FL_2)_\cu$ is a sum of copies of $M$ and $SM$. Therefore, in the
exact sequence
\[
\Hom(L_1,(FL_2)_\cu) \to \Hom(L_1,FL_2) \to
\Hom(L_1,(FL_2)^\cv) \to \Hom(L_1, S(FL_2)_\cu)
\]
the last term vanishes and
\[
\Hom(L_1, (FL_2)^\cv) \iso \Hom(L_1^\cv, (FL_2)^\cv)
\]
identifies with the quotient of $\Hom(L_1, FL_2)$ by the subspace
of morphisms factoring through a sum of copies of $M$ and $SM$.
Since $\Hom(L_1, SM)$ vanishes, this is also the subspace of
morphisms factoring through a sum of copies of $M$. To finish the
proof, it remains to be noticed that under the canonical
equivalence $\cd^b(\ch') \iso \cv$, if $L_2^\cv$ corresponds to
$L_2'$, then the object $(FL_2)^\cv$ does correspond to
$\tau^{-1}_{\ch'} S L'_2$, by lemma~2.14 of \cite{BMR2} or
section~8.1 of \cite{keller}.

\end{proof}

\section{Example on the Kronecker quiver}\label{}
\begin{subsection}{}

We consider the Kronecker quiver $Q$ obtained by choosing the
following orientation of the diagram $\tilde A_1$:
 \begin{equation}
    \xymatrix{1 \ar@/^/[r]^\alpha \ar@/_/[r]_\beta& 2}.
    \end{equation}
As an illustration of theorems \ref{exchange} and \ref{main}, we
give an interpretation of some results of \cite{SZ} in terms of
the cluster category of the Kronecker quiver. We consider
covariant representations. So if $S_1$, $S_2$ are the simple
representations, then
    $$\dim\Hom(S_i,S_i)=1,\;\dim\Ext^1(S_1,S_2)=2,\;\dim\Ext^1(S_2,S_1)=0.$$
Over a field $k$, the (finite-dimensional) indecomposable
representations of the Kronecker quiver are classified as follows:

The postprojective indecomposable modules $U^n$, $n\geq 0$, the
preinjective indecomposable modules $V^n$, $n\geq 0$, and the
family  of indecomposables modules $W^n=W^n(x)$, $n>0$, of the
(regular) tube parametrized by $x\in\P^1(k)$. They are given by
   $$U^n:\;\; \xymatrix{k^n \ar@/^/[r]^\alpha \ar@/_/[r]_\beta& k^{n+1}}\hbox{ with } \alpha=\left[
\begin{array}{ccc}
  I_n   \\
   0
\end{array}
\right],\;\beta=\left[
\begin{array}{ccc}
  0   \\
   I_n
\end{array}
\right]$$

$$V^n:\;\; \xymatrix{k^{n+1} \ar@/^/[r]^\alpha \ar@/_/[r]_\beta& k^n}\hbox{ with } \alpha=\left[
\begin{array}{ccc}
  I_n & 0
\end{array}
\right],\;\beta=\left[
\begin{array}{ccc}
  0   &
   I_n
\end{array}
\right]$$

$$W^n((1\,:\,\lambda)):\;\; \xymatrix{k^n \ar@/^/[r]^\alpha \ar@/_/[r]_\beta& k^n}\hbox{ with }
\alpha=I_n
  ,\;\beta= \left(
     \raisebox{0.5\depth}{
       \xymatrixcolsep{1ex}
       \xymatrixrowsep{1ex}
       \xymatrix{
         \lambda \ar @{-}[ddddrrrr] &
         1 \ar @{-}[dddrrr]
         & 0 \ar @{.}[rr] \ar@{.}[ddrr]& & 0 \ar@{.}[dd]\\
         0 \ar@{.}[ddd] \ar@{.}[dddrrr]\\
         &&&& 0\\
         &&&& 1 \\
         0 \ar@{.}[rrr] & & & 0 & \lambda
       }
     }
   \right),
$$
for $\lambda$ in $k$, and
$$W^n((0\,:\,1)):\;\; \xymatrix{k^n \ar@/^/[r]^\alpha \ar@/_/[r]_\beta& k^n}\hbox{ with } \alpha=
  \left(
     \raisebox{0.5\depth}{
       \xymatrixcolsep{1ex}
       \xymatrixrowsep{1ex}
       \xymatrix{
         0 \ar @{-}[ddddrrrr] &
         1 \ar @{-}[dddrrr]
         & 0 \ar @{.}[rr] \ar@{.}[ddrr]& & 0 \ar@{.}[dd]\\
         0 \ar@{.}[ddd] \ar@{.}[dddrrr]\\
         &&&& 0\\
         &&&& 1 \\
         0 \ar@{.}[rrr] & & & 0 & 0
       }
     }
   \right)
  ,\;\beta= I_n.
$$

Let us calculate the cluster variables of the cluster algebra
$\ca(Q)\subset\Z[x_1^{\pm 1},x_2^{\pm 1}]$ corresponding to the
Kronecker quiver. By theorem \ref{main}, they are given by
$X_{U^n}$, $X_{V^n}$, $n\geq 0$. By duality,   $X_{V^n}$ is
obtained from $X_{U^n}$ by exchanging $x_1$ and $x_2$. So, we just
have to calculate $X_{U^n}$.

 Let $(y_n)_{n\in\N}$ be the sequence given by $y_0=x_2$, $y_1=x_1$,
and $y_{n+2}=X_{U^n}$, $n\geq 0$.  Set $P_1=U^1$, $P_2=U^0$, for
the indecomposable projective modules, then $SP_2\oplus SP_1$, is
the "seed" tilting object and $SP_1\oplus U^0$ is also a tilting
object of the cluster category since the first component of
$\dimv(U^0)$ is zero. Moreover,  for any $n\geq 0$, $U^n\oplus
U^{n+1}$ is easily seen to be a tilting $kQ$-module by applying
recursively the inverse AR-functor to the tilting module
$P_2\oplus P_1$ and the object $SP_1\oplus U^0$ . By applying
theorem \ref{main}, we obtain that  the $y_n$'s are cluster
variables and that
$$\mu_2(\{y_{2n},y_{2n+1}\})=\{y_{2n+2},y_{2n+1}\},\,\,\mu_1(\{y_{2n-1},y_{2n}\})=\{y_{2n+1},y_{2n}\}.$$
In particular, the exchange relations imply that the sequence
$(y_n)_{n\in \N}$ is given by
\begin{equation}\label{induction}
y_0=x_2,\,\,y_1=x_1,\,\, y_{n-1}y_{n+1}=y_n^2+1.
\end{equation}
Note that in the module category $\Ext^1(W^1,P_i)=k$, for $i=1,2$,
which implies
\[
\Ext^1_{\cc_Q}(W^1,P_i)=k \ko
\]
because the $P_i's$ are projective. Applying the
(AR)-autoequivalence $\tau$ in the cluster category, we obtain
$\Ext_{\cc_Q}^1(W^1,U^n)=k$, $n\geq 0$. In the module category, we
have (up to isomorphism) a unique non split exact sequence
$$0\rightarrow U^n\rightarrow U^{n+1}\rightarrow W^1\rightarrow 0.$$
This yields a triangle in the cluster category
$$U^n\rightarrow U^{n+1}\rightarrow W^1\rightarrow SU^n.$$
But, as $SU^n=\tau U^n=U^{n-2}$ in the cluster category, shifting
the triangle gives
$$W^1\rightarrow U^{n-1}\rightarrow U^n\rightarrow SW^ 1.$$
Now, let $w_1:=X_{W^1}=\frac{1+x_1^2+x_2^2}{x_1x_2}$, then theorem
\ref{exchange} implies
$$w_1y_n=y_{n+1}+y_{n-1}.$$
Note that this formula simplifies the initial induction
\ref{induction}. It was obtained in a direct way in \cite{SZ}. We
can calculate the generating series of $(y_n)_{n\in \N}$
 $$\sum_{n\geq 0}y_nt^n=\frac{1-y_{-1}t}{1-w_1t+t^2}\ko$$
where $y_{-1}:=X_{V^0}=\frac{1+x_2^2}{x_1}$.

\end{subsection}

\end{document}